\pdfoutput=1
\RequirePackage{ifpdf}
\ifpdf 
\documentclass[pdftex]{sigma}
\else
\documentclass{sigma}
\fi

\newcommand{\Tr}{\operatorname{Tr}}

\begin{document}

\allowdisplaybreaks

\renewcommand{\thefootnote}{$\star$}

\newcommand{\arXivNumber}{1307.2137}

\renewcommand{\PaperNumber}{040}

\FirstPageHeading

\ShortArticleName{Toda Equations and Piecewise Polynomiality for Mixed Double Hurwitz Numbers}

\ArticleName{Toda Equations and Piecewise Polynomiality\\ for Mixed Double Hurwitz Numbers\footnote{This paper is a~contribution to the Special Issue
on Asymptotics and Universality in Random Matrices, Random Growth Processes, Integrable Systems and Statistical Physics in honor of Percy Deift and Craig Tracy.
The full collection is available at \href{http://www.emis.de/journals/SIGMA/Deift-Tracy.html}{http://www.emis.de/journals/SIGMA/Deift-Tracy.html}}}

\Author{I.P.~GOULDEN~$^\dag$, Mathieu GUAY-PAQUET~$^\ddag$ and Jonathan NOVAK~$^\S$}

\AuthorNameForHeading{I.P.~Goulden, M.~Guay-Paquet and J.~Novak}

\Address{$^\dag$~Department of Combinatorics and Optimization, University of Waterloo,\\
\hphantom{$^\dag$}~200 University Ave.~W., Waterloo, ON, N2L 3G1 Canada}
\EmailD{\href{mailto:ipgoulden@uwaterloo.ca}{ipgoulden@uwaterloo.ca}}
\URLaddressD{\url{https://uwaterloo.ca/math/ian-gouldens-home-page}}

\Address{$^\ddag$~D\'epartement de math\'ematiques, Universit\'e du Qu\'ebec \`a Montr\'eal, \\
\hphantom{$^\ddag$}~C.P.~8888, succ.~Centre-ville, Montr\'eal, Qu\'ebec, H3C 3P8 Canada}
\EmailD{\href{mailto:mathieu.guaypaquet@lacim.ca}{mathieu.guaypaquet@lacim.ca}}

\Address{$^\S$~Department of Mathematics, University of California San Diego,\\
\hphantom{$^\S$}~9500 Gilman Drive, La Jolla, CA 92093-0404 USA}
\EmailD{\href{mailto:jinovak@ucsd.edu}{jinovak@ucsd.edu}}

\ArticleDates{Received February 02, 2016, in f\/inal form April 13, 2016; Published online April 20, 2016}

\Abstract{This article introduces \emph{mixed double Hurwitz numbers}, which interpolate
combinatorially between the classical double Hurwitz numbers studied by Okounkov
and the monotone double Hurwitz numbers introduced recently by Goulden, Guay-Paquet
and Novak. Generalizing a result of Okounkov, we prove that a certain generating series for the
mixed double Hurwitz numbers solves the 2-Toda hierarchy of partial dif\/ferential equations.
We also prove that the mixed double Hurwitz numbers are piecewise polynomial, thereby generalizing
a result of Goulden, Jackson and Vakil.}


\Classification{05A05; 14H70}

\renewcommand{\thefootnote}{\arabic{footnote}}
\setcounter{footnote}{0}

\section{Introduction}
Consider the right Cayley graph of the symmetric group $S(d)$,
as generated by the full conjugacy class of transpositions.
This is a $\binom{d}{2}$-regular graded graph with levels $L_0,L_1,\dots,L_{d-1}$,
where~$L_k$ is the set of permutations with $d-k$ cycles.
Each level $L_k$ decomposes as the disjoint union
of those conjugacy classes in $S(d)$ labelled by Young diagrams with $d-k$ rows.

Let us introduce an edge labelling of the Cayley graph by marking each edge corresponding
to the transposition $\tau=(s\ t)$ with $t$, the larger of the two elements interchanged by
$\tau$. This edge labelling was used by Stanley~\cite{Stanley} and Biane~\cite{Biane} to
study various connections between permutations, parking functions, and noncrossing
partitions\footnote{Stanley and Biane use $s$ as the edge label instead of $t$, but
this is a minor dif\/ference.}.

Given integers $k,l \geq 0$ and partitions $\alpha,\beta \vdash d$,
let $W^{k,l}(\alpha,\beta)$ denote the number of $(k+l)$-step
walks
	\begin{gather*}
		\sigma = \rho\underbrace{(s_1\ t_1) \cdots (s_k\ t_k)}_k \underbrace{(s_{k+1}\ t_{k+1}) \cdots (s_{k+l}\ t_{k+l})}_l
	\end{gather*}		
on the Cayley graph beginning in the conjugacy class $C_\alpha$ and ending in the conjugacy class
$C_\beta$ which satisfy
	\begin{gather*}
		t_1 \leq \dots \leq t_k.
	\end{gather*}	
In words, this monotonicity condition states that
the labels of the edges traversed in the f\/irst $k$ steps of the walk form
a weakly increasing sequence.

While elementary to def\/ine, the numbers $W^{k,l}(\alpha,\beta)$ are related to some
rather sophisticated mathematics. Let $z,t,u,a_1,a_2,\dots,b_1,b_2,\dots$ be
commuting indeterminates, and form the ge\-ne\-rating function
	\begin{gather*}
		\mathbf{W}(z,t,u,A,B) = 1 + \sum_{d=1}^{\infty} \frac{z^d}{d!}
		\sum_{k,l=0}^{\infty} t^k \frac{u^l}{l!} \sum_{\alpha,\beta \vdash d}
		W^{k,l}(\alpha,\beta) p_\alpha(A) p_\beta(B),
	\end{gather*}	
where $p_\alpha(A)$ and $p_\beta(B)$ denote the power-sum symmetric functions
in the variables $A=\{a_1$, $a_2,\dots\}$ and $B=\{b_1,b_2,\dots\}$, respectively. The series
	\begin{gather*}
		\mathbf{H}(z,t,u,A,B) = \log \mathbf{W}(z,t,u,A,B)
	\end{gather*}
is a well-def\/ined element of $\mathbb{Q}[[z,t,u,A,B]]$. Set
	\begin{gather*}
		H^{k,l}(\alpha,\beta) = \bigg{[} z^d t^k \frac{u^l}{l!} p_\alpha(A) p_\beta(B) \bigg{]}\mathbf{H}(z,t,u,A,B),
	\end{gather*}
where $[X]Y$ denotes the coef\/f\/icient of the term $X$ in a series $Y$.

The numbers $H^{0,l}(\alpha,\beta)$ were f\/irst studied by Okounkov \cite{Okounkov},
who called them the \emph{double Hurwitz numbers}. By a classical
construction due to Hurwitz~-- the monodromy construction~-- $H^{0,l}(\alpha,\beta)$ is a weighted count of degree~$d$ branched covers
of the Riemann sphere by a~compact, connected Riemann surface such that the covering map has prof\/ile $\alpha$ over~$0$,
$\beta$~over~$\infty$, and simple ramif\/ication over each of the~$l$th roots of unity. The Riemann--Hurwitz
formula determines the genus of the covering surface in terms of the ramif\/ication data of the covering map:
	\begin{gather*}
		g=\frac{l+2-\ell(\alpha)-\ell(\beta)}{2}.
	\end{gather*}

{\sloppy 	
Verifying and extending a conjecture of Pandharipande~\cite{Pan} in Gromov--Witten theory,
Okoun\-kov proved that the generating function $\mathbf{H}(z,0,u,A,B)$ is a solution
of the $2$-Toda hierarchy of Ueno and Takasaki.
The $2$-Toda hierarchy is a countable collection
of partial dif\/ferential equations, each of which yields a recurrence relation satisf\/ied by the double
Hurwitz numbers. A~construction of the Toda hierarchy may be found in~\cite[Section~4]{OP2}.
Kazarian and Lando~\cite{KL} subsequently showed that, when combined with
the ELSV formula~\cite{ELSV}, Okounkov's result yields a~streamlined proof of the
Kontsevich--Witten theorem relating intersection theory in moduli spaces of curves to
integrable hierarchies.

}

Goulden, Jackson and Vakil~\cite{GJV} gave an alternative interpretation of the
double Hurwitz number $H^{0,l}(\alpha,\beta)$ as counting lattice points in a
certain integral polytope. As a~consequence of this interpretation
and Ehrhart's theorem, it was shown in~\cite{GJV} that, after a simple rescaling,
$H^{0,l}(\alpha,\beta)$ is a piecewise
polynomial function of the parts of $\alpha$ and $\beta$, when $\ell(\alpha)$ and $\ell(\beta)$
are held f\/ixed. Detailed structural properties of this piecewise polynomial behaviour were postulated
in~\cite{GJV}, and subsequently shown to hold by Johnson~\cite{Johnson} using
the combinatorics of the inf\/inite wedge representation of $\mathfrak{gl}(\infty)$.
Shadrin, Spitz and Zvonkine~\cite{SSZ} have recently generalized double Hurwitz numbers
to the setting of the completed cycle theory introduced by Okounkov and Pandharipande~\cite{OP}.
It is shown in~\cite{SSZ} that piecewise polynomiality of double Hurwitz numbers with completed
cycle insertions follows from a suitable modif\/ication of Johnson's arguments.

The \emph{monotone double Hurwitz numbers} $H^{k,0}(\alpha,\beta)$ were introduced by
the present authors in \cite{GGN1}, where it was shown that they are the combinatorial objects
underlying the asymptotic expansion of the Harish--Chandra--Itzykson--Zuber integral, an important special function in random matrix theory.
The paper~\cite{GGN1} was an outgrowth of the works~\cite{MN, Novak:IMRN}, which developed new connections between unitary matrix integrals,
combinatorics, and integrable systems.
Structural properties of the monotone single Hurwitz numbers
$H^{k,0}(\alpha)=H^{k,0}(\alpha,1^d)$
were studied in detail in~\cite{GGN2,GGN3}, where it was shown that they enjoy
a high degree of structural similarity with the classical single Hurwitz numbers $H^{0,l}(\alpha)=H^{0,l}(\alpha,1^d)$.

In this article, we extend the theorems of Okounkov and Goulden--Jackson--Vakil to the more general setting of
the \emph{mixed double Hurwitz numbers} $H^{k,l}(\alpha,\beta)$, which interpolate between the
classical double Hurwitz numbers ($k=0$) and the monotone double Hurwitz numbers~($l=0$).

\begin{theorem}\label{thm:integrability}
The generating function~$\mathbf{H}$ is a solution of the $2$-Toda hierarchy.
\end{theorem}
	
	\begin{theorem}
		\label{thm:polynomiality}
		The mixed double Hurwitz numbers are piecewise
		polynomial.
	\end{theorem}

\section{Toda equations}\label{sec:Toda}

Let us group the transposition generators of $S(d)$ into a matrix,
	\begin{gather*}
		T=\begin{bmatrix}
			(1\ 2) & (1\ 3) & \dots & (1\ d) \\
			{} & (2\ 3) & \dots & (2\ d) \\
			{} & {} & \ddots & \vdots\\
			{} & {} & {} & (d-1\ d)
		\end{bmatrix}.
	\end{gather*}
Denote by
	\begin{gather*}
		J_2 = (1\ 2), \qquad 		J_3 = (1\ 3) + (2\ 3), \qquad
		\dots,\qquad
		J_d = (1\ d) + (2\ d) + \dots + (d-1\ d),
	\end{gather*}		
the column sums of this matrix, viewed as elements of
the group algebra $\mathbb{Q}S(d)$. These elements commute.
Set $J_1:=0$, and introduce the multiset
	\begin{gather*}
		\Xi_d =\{\{J_1,\dots,J_d,0,0,\dots\}\}.
	\end{gather*}

Let $\Lambda$ denote the $\mathbb{Q}$-algebra of symmetric
functions. We will consider the evaluation of the
complete symmetric function $h_\mu \in \Lambda$ indexed by the
$(k,l)$-hook Young diagram $\mu =(k, 1^l)$ on the alphabet $\Xi_d$.
From the def\/inition of the complete symmetric functions, we have
	\begin{align*}
		h_{(k,1^l)}(\Xi_d) &= h_k(\Xi_d) h_1(\Xi_d)^l =\bigg{(} \sum_{2 \leq t_1 \leq \dots \leq t_k \leq d} J_{t_1} \cdots J_{t_k} \bigg{)} \bigg{(} \sum_{t=2}^d J_t \bigg{)}^l \\
		&= \bigg{(}\sum_{2 \leq t_1 \leq \dots \leq t_k \leq d} \bigg{(} \sum_{s_1<t_1} (s_1\ t_1) \bigg{)} \cdots \bigg{(} \sum_{s_k<t_k} (s_k\ t_k) \bigg{)} \bigg{)}
			\bigg{(} \sum_{2 \leq t \leq d} \sum_{s < t} (s\ t) \bigg{)}^l\\
		&= \sum_{\substack{t_1, \dots t_{k+l}=2\\ t_1 \leq \dots \leq t_k\\ s_i < t_i}}^d (s_1\ t_1) \cdots (s_k\ t_k) (s_{k+1}\ t_{k+1}) \cdots (s_{k+l}\ t_{k+l}).
	\end{align*}
Thus
	\begin{gather*}
		W^{k,l}(\alpha,\beta) = [C_{(1^d)}]C_\alpha h_{(k,1^l)}(\Xi_d)C_\beta,
	\end{gather*}	
where we have identif\/ied each conjugacy class in $S(d)$ with the formal sum
of its elements in~$\mathbb{Q}S(d)$. In other words, $W^{k,l}(\alpha,\beta)$
is the normalized character of $C_\alpha h_{(k,1^l)}(\Xi_d)C_\beta$ in the regular representation
of~$\mathbb{Q}S(d)$.

The columns sums of $T$ are known as the \emph{Jucys--Murphy
elements} of $\mathbb{Q}S(d)$.
It was observed by Jucys and Murphy (see \cite{DG} for a proof)
that $e_r(\Xi_d)$, the $r$th elementary symmetric function
evaluated on the alphabet of Jucys--Murphy elements,
is precisely the sum of the permutations
on level $L_r$ of the Cayley graph:
	\begin{gather*}
		e_r(\Xi_d) = \sum_{\substack{\mu \vdash d\\ \ell(\mu) = d-r}} C_\mu.
	\end{gather*}	
In particular, $e_r(\Xi_d)$ belongs to the center $Z(d)$ of $\mathbb{Q}S(d)$.
Since $\Lambda=\mathbb{Q}[e_1,e_2,\dots]$, the substitution
$f \mapsto f(\Xi_d)$ def\/ines a specialization $\Lambda \rightarrow Z(d)$.
In fact, since the levels of the Cayley graph generate~$Z(d)$,
this specialization is surjective~\cite{FH}.

Since $C_\alpha h_{(k,1^l)}(\Xi_d) C_\beta$ belongs to the centre of $\mathbb{Q}S(d)$, we
can calculate its character in the regular representation using the Fourier transform.
Let $(V^\lambda,\rho^\lambda)$, $\lambda \vdash d$, be pairwise non-isomorphic irreducible
representations of $\mathbb{Q}S(d)$, so that the map
	\begin{gather*}
		\sigma \mapsto \big(\rho^\lambda(\sigma)\colon \lambda\vdash d\big)
	\end{gather*}	
def\/ines an algebra isomorphism
	\begin{gather*}
		\mathbb{Q}S(d) \rightarrow \bigoplus_{\lambda \vdash d} \operatorname{End}V^\lambda.
	\end{gather*}
The normalized character of $C_\alpha h_{(k,1^l)}(\Xi_d)C_\beta$ in the regular representation of $\mathbb{Q}S(d)$
is thus
	\begin{align*}
		[C_{(1^d)}] C_\alpha h_{(k,1^l)}(\Xi_d)C_\beta
		&= \sum_{\lambda \vdash d} \Tr \rho^\lambda(C_\alpha h_{(k,1^l)}(\Xi_d) C_\beta) \frac{\dim V^\lambda}{d!} \\
		&= \sum_{\lambda \vdash d} \Tr \rho^\lambda(C_\alpha) \rho^\lambda(h_{(k,1^l)}(\Xi_d)) \rho^\lambda(C_\beta) \frac{\dim V^\lambda}{d!} \\
		&= \sum_{\lambda \vdash d} \Tr \omega^\lambda(C_\alpha) \omega^\lambda(h_{(k,1^l)}(\Xi_d)) \omega^\lambda(C_\beta) \frac{(\dim V^\lambda)^2}{d!},
	\end{align*}	
where, for any $C \in Z(d)$, we denote by $\omega^\lambda(C)$ the unique eigenvalue
of the scalar operator $\rho^\lambda(C) \in \operatorname{End} V^\lambda$, i.e.,
$\rho^\lambda(C) = \omega^\lambda(C)I_{V^\lambda}$. The eigenvalue $\omega^\lambda(C)$ is known
as the \emph{central character} of $C$ in the representation $(V^\lambda,\rho^\lambda)$.

The central character of any conjugacy class $C_\mu$ is given, in terms of the usual character
	\begin{gather*}
		\chi^\lambda_\mu = \Tr \rho^\lambda(\sigma), \qquad \sigma \in C_\mu,
	\end{gather*}	
by the formula
	\begin{gather*}
		\omega^\lambda(C_\mu) = |C_\mu| \frac{\chi^\lambda_\mu}{\dim V^\lambda}.
	\end{gather*}
The central character of any symmetric function $f$ evaluated on $\Xi_d$ is simply
	\begin{gather*}
		\omega^\lambda(f(\Xi_d))=f(\operatorname{Cont}_\lambda),
	\end{gather*}			
the evaluation of $f$ on the multiset of contents of the Young diagram~$\lambda$.
This remarkable result is due to Jucys and Murphy, see~\cite{DG} for a~proof.

Recalling that the Schur functions have the expansion
	\begin{gather*}
		s_\lambda = \ \sum_{\mu \vdash d} \frac{|C_\mu|}{d!}\chi^\lambda_\mu p_\mu,
	\end{gather*}
where $\lambda \vdash d$, the generating function $\mathbf{W}= \mathbf{W}(z,t,u,A,B)$ may be
rewritten as follows:
	\begin{align*}
\mathbf{W}& = 1+\sum_{d=1}^{\infty}\frac{z^d}{d!} \sum_{k,l=0}^{\infty} t^k \frac{u^l}{l!}
		\sum_{\alpha,\beta,\lambda \vdash d} W^{k,l}(\alpha,\beta) p_\alpha(A) p_\beta(B) \\
		&= 1+\sum_{d=1}^{\infty}\frac{z^d}{d!} \sum_{k,l=0}^{\infty} t^k \frac{u^l}{l!}
		\sum_{\alpha,\beta,\lambda \vdash d} \bigg{(} \frac{1}{d!}\sum_{\lambda \vdash d}
		|C_\alpha| \chi^\lambda_\alpha h_{(k,1^l)}(\operatorname{Cont}_\lambda) |C_\beta| \chi^\lambda_\beta\bigg{)} p_\alpha(A) p_\beta(B) \\
		&= 1+\sum_{d=1}^{\infty}z^d \sum_{k,l=0}^{\infty} t^k \frac{u^l}{l!} \sum_{\lambda \vdash d} h_{(k,1^l)}(\operatorname{Cont}_\lambda)
		\bigg{(} \sum_{\alpha \vdash d} \frac{|C_\alpha|}{d!} \chi^\lambda_\alpha p_\alpha(A)\bigg{)}
		\bigg{(} \sum_{\beta \vdash d} \frac{|C_\beta|}{d!} \chi^\lambda_\beta p_\beta(B)\bigg{)} \\
		&=1+\sum_{d=1}^{\infty}z^d \sum_{k,l=0}^{\infty} t^k \frac{u^l}{l!} \sum_{\lambda \vdash d} h_{(k,1^l)}(\operatorname{Cont}_\lambda)s_\lambda(A)s_\lambda(B) \\
		&=1+\sum_{d=1}^{\infty}z^d \sum_{\lambda \vdash d} s_\lambda(A)s_\lambda(B) \bigg{(} \sum_{k=0}^{\infty} t^k h_k(\operatorname{Cont}_\lambda) \bigg{)}
			\bigg{(} \sum_{l=0}^{\infty} \frac{u^l}{l!} h_1(\operatorname{Cont}_\lambda)^l \bigg{)} \\
		&= \sum_{\lambda \in \mathcal{Y}} Y(\lambda) s_\lambda(A)s_\lambda(B).
	\end{align*}	
Here $\mathcal{Y}$ is the set of all Young diagrams (including the empty diagram) and
	\begin{gather*}
		Y(\lambda) = Y(\lambda;z,t,u) = \prod_{\Box \in \lambda} \frac{ze^{c(\Box)u}}{1-c(\Box)t},
	\end{gather*}
where for any Young diagram $\lambda$ and cell $\Box$ in $\lambda$, $c(\Box)$ denotes the content of this
cell, i.e., its column index less its row index. By convention, an empty product equals $1$.

In order to complete the proof of Theorem~\ref{thm:integrability}, we appeal to
the following result of Orlov and Shcherbin~\cite{OS}, and Carrell~\cite{Carrell}.

	\begin{theorem}
		\label{thm:contentProduct}
		Let $\{y_k \colon k \in \mathbb{Z}\}$ be a set of variables indexed by the integers, and
		set
			\begin{gather*}
				Y(\lambda) = \prod_{\Box \in \lambda} y_{c(\Box)}
			\end{gather*}
		for each $\lambda \in \mathcal{Y}$. The series
			\begin{gather*}
				\log \bigg{(} \sum_{\lambda \in \mathcal{Y}}
				Y(\lambda) s_\lambda(A) s_\lambda(B)\bigg{)}
			\end{gather*}
		is a solution of the $2$-Toda lattice hierarchy in the variables
		$p_1(A),p_2(A),\dots$ and $p_1(B)$,\linebreak $p_2(B),\dots$.
	\end{theorem}

Solutions of the Toda equations of the form described in Theorem \ref{thm:contentProduct}
are known as \emph{diagonal content-product solutions}.
Our computations above show that the generating function
	\begin{gather*}
		\mathbf{H}(z,t,u,A,B) = \log \mathbf{W}(z,t,u,A,B)
	\end{gather*}
of the mixed double Hurwitz numbers is a diagonal content-product
solution, with
	\begin{gather*}
		y_k = y_k(z,t,u) = \frac{ze^{ku}}{1-kt}.
	\end{gather*}	
In particular, Okounkov's generating function $\mathbf{H}(z,0,u,A,B)$ for the classical double Hurwitz
numbers is the diagonal content product solution with
	\begin{gather*}
		y_k = ze^{ku},
	\end{gather*}	
while the generating function $\mathbf{H}(z,t,0,A,B)$ for the monotone double Hurwitz numbers
is the diagonal content product solution with
	\begin{gather*}
		y_k = \frac{z}{1-kt}.
	\end{gather*}

\section{Piecewise polynomiality}
\label{sec:Polynomiality}
To prove Theorem \ref{thm:polynomiality}, let us begin by formulating precisely the
relationship between the num\-bers~$W^{k,l}(\alpha,\beta)$ and $H^{k,l}(\alpha,\beta)$.
To this end, consider the generating functions~$\mathbf{W}$ and~$\mathbf{H}$ as
elements of~$\mathbb{Q}[[t,u,A,B]][[z]]$ by writing
	\begin{gather*}
		\mathbf{W}(z) = 1+ \sum_{d=1}^{\infty} \frac{z^d}{d!} \mathbf{W}_d(t,u,A,B), \qquad
		\mathbf{H}(z) = \sum_{d=1}^{\infty} \frac{z^d}{d!} \mathbf{H}_d(t,u,A,B),
	\end{gather*}
where, for each $d \geq 1$,
	\begin{gather*}
		 \mathbf{W}_d(t,u,A,B) = \sum_{k,l=0}^{\infty} t^k\frac{u^l}{l!} \sum_{\alpha,\beta \vdash d} W^{k,l}(\alpha,\beta) p_\alpha(A)p_\beta(B), \\
		 \mathbf{H}_d(t,u,A,B) = d!\sum_{k,l=0}^{\infty} t^k\frac{u^l}{l!} \sum_{\alpha,\beta \vdash d} H^{k,l}(\alpha,\beta) p_\alpha(A)p_\beta(B).
	\end{gather*}	
Then, by the exponential formula, we have
	\begin{gather*}
		 \mathbf{W}_d(t,u,A,B) = \sum_{r=1}^d\ \sum_{P_1 \sqcup \dots \sqcup P_r}\ \prod_{j=1}^r |P_j|! \mathbf{H}_{|P_j|}(t,u,A,B)
	\end{gather*}		
for each $d \geq 1$, as an identity in $\mathbb{Q}[[t,u,A,B]]$.
In this identity, the inner sum is over partitions $P_1 \sqcup \dots \sqcup P_r$ of $\{1,\dots,d\}$
into $r$ disjoint nonempty sets. Equating the coef\/f\/icient of $t^k\frac{u^l}{l!}$ on either side
of this identity, we obtain
	\begin{align*}
		\sum_{\alpha,\beta \vdash d} W^{k,l}(\alpha,\beta)p_\alpha(A)p_\beta(B) =
		\sum_{r=1}^d\ \!\! \sum_{P_1 \sqcup \dots \sqcup P_r}\ \!\! \prod_{j=1}^r |P_j|!
				\!\!\!\! \sum_{\zeta^j,\eta^j \vdash |P_j|} \!\!\!\! H^{k,l}(\zeta^j,\eta^j) p_{\eta^j}(A) p_{\zeta^j}(B),
	\end{align*}	
for each $d \geq 1$ and $k,l \geq 0$, as an identity in $\mathbb{Q}[[A,B]]$.
Note that the sum on the right hand side depends only on the overall block
structure of the set partition $P_1 \sqcup \dots \sqcup P_r$, and not on the internal structure of the individual blocks.
We can thus replace the right hand side by a sum over integer partitions,
	\begin{gather*}
		 \sum_{\alpha,\beta \vdash d} W^{k,l}(\alpha,\beta)p_\alpha(A)p_\beta(B)\\
		 \qquad {}= \sum_{\theta \vdash d} c_\theta \sum_{\substack{(\zeta^1,\dots,\zeta^{\ell(\theta)})\\ \zeta^j \vdash \theta_j}}
		\sum_{\substack{(\eta^1,\dots,\eta^{\ell(\theta)})\\ \eta^j \vdash \theta_j}} \bigg{(} \prod_{j=1}^{\ell(\theta)} H^{k,l}(\zeta^j,\eta^j) \bigg{)}
		p_{\zeta^1 \cup \dots \cup \zeta^{\ell(\theta)}}(A) p_{\eta^1 \cup \dots \cup \eta^{\ell(\theta)}}(B),
	\end{gather*}	
where the coef\/f\/icient $c_\theta$ is given by
	\begin{gather*}
		c_\theta = |f^{-1}(\theta)| \prod_{i=1}^{\ell(\theta)} \theta_i!
	\end{gather*}	
and $f$ denotes the surjection
	\begin{gather*}
		P_1 \sqcup \dots \sqcup P_r \mapsto \Big( \max_{1 \leq j \leq r} |P_j|, \dots, \min_{1 \leq j \leq r} |P_j|\Big)
	\end{gather*}	
from set partitions of $\{1,\dots,d\}$ onto integer partitions of $d$. The internal sums run over
sequences of partitions whose $j$th element is a partition of the $j$th part of~$\theta$, and if $(\mu^1,\dots,\mu^k)$ is any
sequence of partitions then $\mu^1 \cup \dots \cup \mu^k$ denotes the partition obtained by arranging the parts of
$\mu^1,\dots,\mu^k$ in weakly decreasing order.

Extracting the coef\/f\/icient of $p_\alpha(A)p_\beta(B)$ on
each side of the above identity in $\mathbb{Q}[[A,B]]$, we obtain the numerical identity
	\begin{gather*}
		W^{k,l}(\alpha,\beta) = \sum_{\theta \vdash d} c_\theta
		 \sum_{\substack{(\zeta^1,\dots,\zeta^{\ell(\theta)})\\ \zeta^j \vdash \theta_j \\ \zeta^1 \cup \dots \cup \zeta^{\ell(\theta)}=\alpha}}\
		 \sum_{\substack{(\eta^1,\dots,\eta^{\ell(\theta)})\\ \eta^j \vdash \theta_j \\ \eta^1 \cup \dots \cup \eta^{\ell(\theta)}=\beta}}
		 \prod_{j=1}^{\ell(\theta)} H^{k,l}(\zeta^j,\eta^j).
	\end{gather*}
The sum on the right hand side receives non-zero contributions from Young diagrams $\theta$ whose rows
can be obtained by gluing together rows of $\alpha$, and by gluing together rows of $\beta$. This is
possible if and only if, for each row $\theta_k$ of $\theta$, there exist subsets $I_k,J_k \subseteq \{1,\dots,d\}$ such
that
	\begin{gather}
		\label{eqn:RowConstraint}
		\sum_{i \in I_k} \alpha_i = \sum_{j \in J_k} \beta_j = \theta_k.
	\end{gather}	
In the case of the one-row Young diagram $\theta=(d)$, the required
sets are simply $I_1=J_1=\{1,\dots,d\}$, so we have
	\begin{gather*}
		W^{k,l}(\alpha,\beta) = d!H^{k,l}(\alpha,\beta) + \cdots,
	\end{gather*}	
where the ellipsis stands for contributions from Young diagrams $\theta$ that have
at least two rows. These contributions are zero unless the constraint~\eqref{eqn:RowConstraint} is met for each row of $\theta$.
We thus conclude that $W^{k,l}(\alpha,\beta)$ and
$H^{k,l}(\alpha,\beta)$ typically agree, up to a~factor of~$d!$.

We can give a geometric interpretation of the above as follows.
Fix two positive integers $m$ and $n$ and consider the
convex region $\mathfrak{R}_{m,n}$ given by
	\begin{gather*}
	\bigg\{(x_1,\dots,x_m,y_1,\dots,y_n)\colon x_1 \geq \dots \geq x_m >0,\ y_1 \geq \dots \geq y_n >0,\ \sum_{i=1}^m x_i = \sum_{j=1}^n y_j\bigg\}
	\end{gather*}	
in Euclidean space $\mathbb{R}^{m+n}$.
Pairs of partitions $(\alpha,\beta)$ such that
	\begin{gather*}
		|\alpha|=|\beta|,\qquad \ell(\alpha)=m, \qquad \ell(\beta)=n
	\end{gather*}	
may be identif\/ied with lattice points in $\mathfrak{R}_{m,n}$ via
	\begin{gather*}
		(\alpha,\beta) \mapsto (\alpha_1,\dots,\alpha_m,\beta_1,\dots,\beta_n).
	\end{gather*}
Given proper nonempty subsets $I \subset \{1,\dots,m\}$ and $J \subset \{1,\dots,n\}$,
def\/ine a hyperplane $\mathfrak{W}_{IJ}$ in $\mathbb{R}^{m+n}$ by
	\begin{gather*}
		\mathfrak{W}_{IJ} = \bigg\{(x_1,\dots,x_m,y_1,\dots,y_n) \in \mathbb{R}^{m+n}\colon \sum_{i \in I} x_i = \sum_{j \in J} y_j\bigg\}.
	\end{gather*}	
The hyperplanes $\mathfrak{W}_{IJ}$, as $I$ ranges over proper nonempty subsets of $\{1,\dots,m\}$ and~$J$ ranges over proper nonempty subsets of $\{1,\dots,n\}$, constitute the \emph{resonance arrangement}
of~\cite{Johnson}. A~chamber~$\mathfrak{c}$ of the resonance arrangement is a connected component of the
complement of a~hyperplane~$\mathfrak{W}_{IJ}$ in~$\mathfrak{R}_{m,n}$. On any chamber
$\mathfrak{c}$ of the resonance arrangement, we have that
	\begin{gather*}
		W^{k,l}(\alpha,\beta) = d! H^{k,l}(\alpha,\beta)
	\end{gather*}	
for all lattice points $(\alpha,\beta) \in \mathfrak{c}$.

Theorem~\ref{thm:polynomiality} claims that for each chamber $\mathfrak{c}$
there exists a polynomial $p_{\mathfrak{c}}^{k,l}$ in~$m+n$ variables such that
	\begin{gather*}
		H^{k,l}(\alpha,\beta) = p^{k,l}_{\mathfrak{c}}(\alpha_1,\dots,\alpha_m,\beta_1,\dots,\beta_n)
	\end{gather*}
for all lattice points $(\alpha,\beta) \in \mathfrak{c}$. From the above discussion and Section~\ref{sec:Toda},
we know that $H^{k,l}(\alpha,\beta)$ is given by the character formula
	\begin{gather}
		\label{eqn:CharacterFormula}
		H^{k,l}(\alpha,\beta) = \frac{|C_\alpha|}{d!} \frac{|C_\beta|}{d!} \sum_{\lambda \vdash d} \chi^\lambda_\alpha
		h_{(k,1^l)}(\operatorname{Cont}_\lambda) \chi^\lambda_\beta, \qquad d=|\alpha|=|\beta|,
	\end{gather}	
on any chamber of the resonance arrangement. In order to deduce Theorem~\ref{thm:polynomiality} from this character formula, we appeal to a~recent result of
Shadrin, Spitz and Zvonkine~\cite{SSZ} which asserts the piecewise polynomiality of
a general class of sums of the above form.

Let $\mathbb{C}^{\mathcal{Y}}$ denote the algebra of all functions $\mathcal{Y} \rightarrow \mathbb{C}$.
Following Olshanski~\cite[Proposition~2.4]{Olshanski},
we def\/ine the algebra $\mathbb{A}$ of \emph{regular functions} on Young diagrams to
be the subalgebra of~$\mathbb{C}^{\mathcal{Y}}$ generated by the functions
	\begin{gather*}
		\lambda \mapsto f(\operatorname{Cont}_\lambda), \qquad f \in \Lambda,
	\end{gather*}
together with the function $\lambda \mapsto |\lambda|$. Consider the transform
	\begin{gather*}
		S\colon \ \mathbb{C}^{\mathcal{Y}} \rightarrow \mathbb{C}^{\mathcal{Y} \times \mathcal{Y}}
	\end{gather*}	
from functions on Young diagrams to functions on pairs of Young diagrams def\/ined by
	\begin{gather*}
		S^f(\alpha,\beta) = \frac{|C_\alpha|}{d!} \frac{|C_\beta|}{d!} \sum_{\lambda \vdash d} \chi^\lambda_\alpha
		f(\lambda) \chi^\lambda_\beta,
	\end{gather*}
where $d=|\alpha|=|\beta|$. This def\/inition assumes that $\alpha$, $\beta$ have the same size; if $|\alpha|>|\beta|$ or
vice versa, complete the smaller diagram by adding unicellular rows.
As detailed by Olshanski~\cite{Olshanski},
the algebra $\mathbb{A}$ of regular functions on Young diagrams may equivalently
be described as the algebra of functions on~$\mathcal{Y}$ such that~$f(\lambda)$
is a \emph{shifted symmetric function} of the row lengths $\lambda_1,\lambda_2,\dots$
of~$\lambda$. Moreover, the shifted power-sums
	\begin{gather*}
		p_k^*(\lambda) = \sum_{i \geq 1} \left[ \left(\lambda_i-i+ \frac{1}{2}\right)^k - \left(-i+\frac{1}{2}\right)^k\right]
	\end{gather*}	
generate $\mathbb{A}$ as a polynomial ring, i.e., the functions
	\begin{gather*}
		p_\mu^*(\lambda) = \prod_{i=1}^{\ell(\mu)} p_{\mu_i}^*(\lambda)
	\end{gather*}
form a linear basis of $\mathbb{A}$ as $\mu$ ranges over $\mathcal{Y}$.
Building on work of Johnson~\cite{Johnson}, Shadrin, Spitz, and Zvonkine~\cite{SSZ} have used the semi-inf\/inite wedge space formalism to demonstrate
that the $S$-transform of $p_\mu^*$ is piecewise polynomial for each $\mu \in \mathcal{Y}$
(see Theorem~6.3 and Remark~6.5
in~\cite{SSZ}). Thus, the piecewise polynomiality of the $S$-transform
of an arbitrary regular function is a~direct consequence
of the results of~\cite{SSZ}:

	\begin{theorem}
		\label{thm:SSZ}
		Given a regular function
		$f \in \mathbb{A}$, positive integers~$m$ and~$n$, and a~chamber $\mathfrak{c}$
		of the resonance arrangement in $\mathfrak{R}_{m,n}$, there
		exists a polynomial $p^f_{\mathfrak{c}}$
		in $m+n$ variables such that
			\begin{gather*}
				S^f(\alpha,\beta) =
				p^f_{\mathfrak{c}}(\alpha_1,\dots,\alpha_m,\beta_1,\dots,\beta_n)
			\end{gather*}	
		\noindent
		for all lattice points $(\alpha,\beta) \in \mathfrak{c}$.
	\end{theorem}

From the character formula~\eqref{eqn:CharacterFormula},
we see that Theorem \ref{thm:polynomiality} follows by applying
Theorem~\ref{thm:SSZ} with~$f$ the regular function
	\begin{gather*}
		\lambda \mapsto h_{(k,1^l)}(\operatorname{Cont}_\lambda).
	\end{gather*}
Let us close by remarking that structural properties of the polynomials
representing mixed double Hurwitz numbers can probably be obtained
by applying the methods of~\cite{Johnson} and~\cite{SSZ}. In view of recent
progress in monotone Hurwitz theory~\cite{DDM}, this seems to be a~very
interesting topic for future research.

\pdfbookmark[1]{References}{ref}
\LastPageEnding

\end{document}